\documentclass[12,reqno]{amsart}%
\usepackage{amsfonts}
\usepackage{amsthm}
\usepackage[centertags]{amsmath}
\usepackage{amssymb}
\usepackage{enumerate}
\usepackage{graphicx}
\usepackage[bookmarks=true, pdfstartview={FitH}, pdftitle={Paper},
pdfauthor={Sadia Arshad}]{hyperref}%
\providecommand{\U}[1]{\protect\rule{.1in}{.1in}}
\theoremstyle{plain}
\newtheorem{theorem}{Theorem}[section]

\newtheorem{example}[theorem]{Example}

\newtheorem{remark}[theorem]{Remark}

\numberwithin{equation}{section}
\begin{document}
\title[Dichotomy of Poincare Map]{Dichotomy of Poincare Map and boundedness of solutions of certain non-autonomous periodic Cauchy problems}
\author {Akbar Zada, Sadia Arshad, Gul Rahmat and Aftab Khan}
\address{Department of Mathematics, University of Peshawar, Peshawar, Pakistan}
\email{akbarzada@upesh.edu.pk}
\address{Government College University, Abdus Salam School of Mathematical Sciences,
(ASSMS), Lahore, Pakistan}
\email{sadia$\_$735@yahoo.com}
\address{Government College University, Abdus Salam School of Mathematical Sciences,
(ASSMS), Lahore, Pakistan}
\email{gulassms@gmail.com}
\address{Government College University, Abdus Salam School of Mathematical Sciences,
(ASSMS), Lahore, Pakistan}
\email{aftabm84@gmail.com}

\date{}
\dedicatory{}
\begin{abstract}
In this paper we study the dichotomy of the Poincar$\acute{e}$ map and
give the relations between the dichotomy of the Poincar$\acute{e}$ map
and boundedness of solutions of the following periodic Cauchy problems
$$
\left\{
  \begin{array}{ll}
    \dot{X}(t)=A(t)X(t)+e^{i\mu t}Pb, \quad t\ge 0\\
X(0)=0,
  \end{array}
\right.$$
and
$$
\left\{
  \begin{array}{ll}
    \dot{X}(t)=-X(t)A(t)+e^{i\mu t}(I-P)b, \quad t\ge 0\\
X(0)=0,
  \end{array}
\right.$$
where $A(t)$ is a square size matrix of order $m$, $\mu$ is any real number, $b$ is a non zero vector in $\mathbb{C}^{m}$ and $P$ is an orthogonal projection.
\end{abstract}
\maketitle
\section{introduction}
The aim of this paper is to study the relation between dichotomy of Poincar$\acute{e}$ map and boundedness of the solutions of the $q$-periodic $(q>0)$ Cauchy problems in the continuous case. For a well-posed non-autonomous Cauchy problem
$$
\left\{
  \begin{array}{ll}
    \dot{x}(t)=A(t)x(t)+e^{i\mu t}I, \quad t\ge 0\\
x(0)=0,
  \end{array}
\right.\eqno{(A(t), \mu,I, 0)}$$
where $A(t)$ an $m\times m$ matrix, the solution leads to an evolution family ${\mathcal{U}}=\{U(t,s), t\geq s\geq 0\},$ i.e. $U(t,s)U(s,r)=U(t,r)$ and $U(t,t)=I$ for all $t\geq s\geq r\geq 0.$
When the Cauchy problem $(A(t), \mu, Pb, 0)$ is $q$-periodic, i.e. $A(t+q)=A(t)$ for all $t\geq 0$, then the family $\mathcal{U}$ is $q$-periodic as well, i.e. $U(t+q,s+q)=U(t,s)$ for all $t\geq s\geq 0$.
It is given in \cite{[SB]} that the evolution family ${\mathcal{U}}$ is uniformly exponentially stable if and only if
the spectral radius of $U(q, 0)$ is less than one,  i.e.
$$r(U(q, 0)):=\sup\{|\lambda|,\,\,\lambda\in\sigma(U(q, 0))\}=\inf\limits_{n\geq 1}\|U(q, 0)^{n}\|^{\frac{1}{n}}<1.$$
 We show that $U(q,0)$ is dichotomic if for each $\mu \in \mathbb{R}$ the matrices
  $$\Phi_{\mu}(q)=\int_{0}^{q} U(q,s)e^{i\mu s}ds\,\,\,\, \text{      and      }\,\,\,\, \Psi_{\mu}(q)=\int_{0}^{q} U^{-1}(q,s)e^{i\mu s}ds$$
  are invertible and there exits a projection $P$ which commutes with $U(q,0),\,\Phi_{\mu}(q)$ and $\Psi_{\mu}(q)$ such that for each real $\mu\in \mathbb{R}$ and each vector $b\in \mathbb{C}^{m},$ the solutions of the Cauchy Problems $(A(t),\mu,Pb,0)$ and $(-A(t), \mu, (I-P)b, 0)$ are bounded on $\mathbb{R}_{+}.$ We give an example that invertibility of the matrices $\Phi_{\mu}(q)$ and $\Psi_{\mu}(q)$ is necessary condition and boundedness of the Cauchy problems $(A(t),\mu,Pb,0)$ and $(-A(t), \mu, (I-P)b, 0)$ is not sufficient for the dichotomy of $U(q,0).$

In \cite{[SB]} and \cite{[ABNZ11]} stability of the poincr$\acute{e}$ map have been studied in the discrete and continuous case respectively. These papers give a connection between stability of the Poincar$\acute{e}$ map and boundedness of the solutions of Cauchy problems. Results regarding the dichotomy of a matrix have been  discussed in \cite{[BZ]} and \cite{[Z1]}. For connection between stability and periodic systems see the papers \cite{[SB]}, \cite{[ABNZ11]}, \cite{[BZ09]} and \cite{[BCDS]}. General theory of dichotomy of infinite dimensional systems has given in the monograph \cite{[CY]}.

The paper is organized as follows: In section 2 we recall basic well known properties of the evolution family. In section 3 we established the results regarding the connection between dichotomy of the Poincar$\acute{e}$ map $U(q,0)$ and boundedness of solutions for some periodic Cauchy problems.
\section{Preliminary Results}
Let $X$ be a Banach space and let $\mathcal{L}(X)$ be the space of all bounded
linear operators acting on $X$. The norm in $X$ and in $\mathcal{L}(X)$ is denoted
by the same symbol $\|.\|$.

A family $\mathcal{U} = \{U(t, s) : t \geq s\geq0 \} \subseteq \mathcal{L}(X)$ is called evolution family if the following properties are satisfied\\
$(i)$ $U(t, t) = I,\text{  for all  } t\in \mathbb{R}_{+},$\\
$(i)$ $ U(t, s)U(s, r) = U(t, r) \text{  for all  } t \geq s\geq r\geq0,$\\
where $I$ denote the identity operator on $\mathcal{L}(X)$. If the later condition is satisfied for all $t,\,s,\,r\in \mathbb{R}_{+}$ then we say that $\mathcal{U}$ is reversible evolution family on $X.$ In this case $U(t,s)$ is invertible for all $t,\,s\in \mathbb{R}_{+}.$
An evolution family $\mathcal{U}$ is called strongly continuous if for each $x\in X$ the map
$$(t, s) \to U(t, s)x : (t, s) \in \mathbb{R}^{2}  \to X $$
is continuous for all $t \geq s \geq 0$. Such a family is called $q$-periodic (with some $q > 0$) if
$$U(t + q, s + q) = U(t, s), \text {   for all }t \geq s \geq 0.$$
Clearly, a $q$-periodic evolution family also satisfies\\
$(i)$ $U(pq + v, pq + u) = U(v, u), \text{ for all  } p \in \mathbb{N}, \text{ for all  } v\geq u \geq 0,$\\
$(ii)$ $U(pq, rq) = U((p - r)q, 0) = U(q, 0)^{p-r}, \text{ for all  } p,r \in \mathbb{N},\, p\geq r.$

Let $\{U(t, s) : t \geq s \geq0\}$ be $q$-periodic evolution family then the operator $U(q,0)$ is called Poincar$\acute{e}$ map or monodromy operator.

The family $\mathcal{U}$ is called uniformly exponentially stable if there exist two positive constants $N$ and $\omega$ such that
$$\|U(t, s)\| \leq Ne^{-\omega(t-s)}, \text{  for all } t \geq s \geq 0.$$

The set of all $m\times m$ matrices having complex entries would be denoted by ${\mathcal{M}}(m,\mathbb{C})$. Assume that the map $t\mapsto A(t):\mathbb{R}\mapsto {\mathcal{M}}(m,\mathbb{C})$ is continuous. Then the Cauchy Problem
$$
\left\{
  \begin{array}{ll}
    \dot{X}(t)=A(t)X(t),\quad t\in \mathbb{R} \\
X(0)=I,
  \end{array}
\right.\eqno{(1)}
$$
has a unique solution denoted by $\Phi(t).$ It is well known that $\Phi(t)$ is an invertible matrix and that its inverse is the unique solution of the Cauchy Problem
$$
\left\{
  \begin{array}{ll}
    \dot{X}(t)=-X(t)A(t),\quad t\in \mathbb{R} \\
X(0)=I.
  \end{array}
\right.\eqno{(2)}
$$
Set $U(t,s):=\Phi(t)\Phi^{-1}(s)\text{ for all }t,s\in \mathbb{R}.$ Obviously, the family ${\mathcal{U}}=\{U(t,s), t,s\in\mathbb{R}\},$ has the following properties:\\

(i) $U(t,t)=I, \text{ for all }t\in \mathbb{R};$\\

(ii) $U(t,s)=U(t,r)U(r,s)\text{ for all }t, s, r\in \mathbb{R};$\\

(iii) $U(t,s)$ is invertible $\text{ for all }t,s\in \mathbb{R};$\\

(iv) $\frac{\partial}{\partial t}U(t,s)=A(t)U(t,s)\text{ for all }t,s\in \mathbb{R};$\\

(v) $\frac{\partial}{\partial s}U(t,s)=-U(t,s)A(s)\text{ for all }t,s\in \mathbb{R};$\\

(vi) The map $(t,s)\mapsto U(t,s):\mathbb{R}^{2}\to \mathcal{M}(m,\mathbb{C})\text{ is continuous }.$\\

If, in addition, the map $A(\cdot)$ is $q$- periodic, for some positive number $q,$ then:\\

(vii) $U(t+q,s+q)= U(t,s)\text{ for all   }t,s\in \mathbb{R};$\\

(viii) $\text{ There exists  } \omega\in \mathbb{R} \text{ and  }  M\geq1 \text{  such that }$
$$\|U(t,s)\|\leq  Me^{\omega(t-s)}, \text{  for all  } \quad t\geq s,$$\\
i.e. the family $\mathcal{U}$ has an exponential growth.

For a given real number $\mu$ and a given family $(A(t))$  we consider the Cauchy Problem
$$
\left\{
  \begin{array}{ll}
    \dot{X}(t)=A(t)X(t)+e^{i\mu t}I, \quad t\ge 0\\
X(0)=0,
  \end{array}
\right.\eqno{(A(t), \mu, I, 0)}$$
and the differential matrix system $$ \dot{X}(t)=A(t)X(t),\quad t\in\mathbb{R}.\eqno{(A(t))}$$

Obviously, the solution of $(A(t),\mu,I,0)$ is given by
$$\Phi_{\mu}(t)=\int_{0}^{t} U(t,s)e^{i\mu s}ds.$$
Now we define
$$V(t,s):=U^{-1}(t,s)=\Phi(s)\Phi^{-1}(t),\,\,t,s\in\mathbb{R} $$
then the family ${\mathcal{V}}=\{V(t,s), t,s\in\mathbb{R}\}$ is an evolution family if
\begin{equation}\label{3}
\Phi(t)\Phi^{-1}(s)=\Phi^{-1}(s)\Phi(t) \text{  for all } t,s\in\mathbb{R}.
\end{equation}
Throughout the paper we assume that equation (\ref{3}) is satisfied for all $t,s\in\mathbb{R}$.\\
Consider the Cauchy problem
$$
\left\{
  \begin{array}{ll}
    \dot{Y}(t)=-Y(t)A(t)+e^{i\mu t}I, \quad t\ge 0\\
Y(0)=0.
  \end{array}
\right.\eqno{(-A(t), \mu, I, 0)}$$

The solution of $(-A(t),\mu,I,0)$ is given by
$$\Psi_{\mu}(t)=\int_{0}^{t} V(t,s)e^{i\mu s}ds.$$

Let $p_L$ be the characteristic polynomial associated to the
matrix $L\in \mathcal{M}(m,\mathbb{C})$ and let
$\sigma{(L)}= {\{\lambda_1,\lambda_2,\dots ,\lambda_k \}} $,
$k\leq {m}$ be its spectrum. There exist integer numbers  $
m_1, m_2, \dots , m_k \geq{1}$ such that
\[
p_L{(\lambda)} = (\lambda-\lambda_1)^{m_1}
(\lambda-\lambda_2)^{m_2} \dots (\lambda-\lambda_k)^{m_k} ,\quad
m_1 + m_2 + \dots + m_k = m .
\]
Let $j\in \{1,2,\dots,k\}$ and
$Y_{j}:= \ker(L-\lambda_{j}I)^{m_{j}}$ then in \cite{[BZ]}\label{2.2} we have the following important theorem which is useful latter on.
\begin{theorem}\label{2.100}  {\it
For each $z\in \mathbb {C}^m$ there exists $y_{j}\in Y_{j}$,  $j \in {\{1,2,\dots,k\}}$  such that
\[
L^{n}z= L^{n}y_1 + L^{n}y_2 + \dots +  L^{n}y_k.
\]
Moreover, if  $  y_{j}(n):= L^{n}y_{j}$  then $y_{j}(n)\in
Y_{j}$ for all $ n\in\mathbb{Z_+} $ and there exist a
$\mathbb{C}^m$-valued polynomials  $p_{j}(n)$  with
$\deg{(p_{j})} \leq {m_{j}-1}$ such that
\begin{equation*}
  y_{j}(n)= \lambda_{j}^{n}p_{j}(n),\quad n\in \mathbb{Z_+} , \ j \in  \{1,2,\dots,k\}.
\end{equation*}}
\end{theorem}

Indeed from the Hamilton-Cayley theorem and using the well known fact that
$$
\ker[pr(L)]= \ker[p(L)]\oplus \ker[r(L)]
$$
whenever the complex valued polynomials $p$ and $r$ are relatively prime,
 follows
\begin{equation*}
\mathbb{C}^m=Y_1\oplus Y_2\oplus \dots \oplus Y_k.\eqno(3.1)
\end{equation*}
Let $z\in{\mathbb{C}^m} $. For each $j\in {\{1,2,\dots,k\}}$ there
exists a unique  $y_{j}\in Y_{j}$ such that
\[
z= y_1 + y_2 + \dots + y_k
\]
and then
\begin{equation*}
L^{n}z= L^{n}y_1 +  L^{n}y_2 + \dots +  L^{n}y_k,
\quad  n \in \mathbb{Z_+}.
\end{equation*}
\section{Dichotomy and Boundedness}
Let us denote $\Gamma_{1}
= {\{ z \in\mathbb{C} : |z|=1\}}$, $\Gamma_1^{+}:={ \{z \in\mathbb{C} :\mathop |z| > 1 \}}$ and
$\Gamma_{1}^{-}: = {\{z \in\mathbb{C} :\mathop |z| < 1\}}$. Clearly
$\mathbb{C}=\Gamma_1 \cup \Gamma_1^{+} \cup \Gamma_1^{-}$.

A matrix $L$ is called:
\begin{itemize}
\item[(i)]  \textit{stable} if $\sigma(L)$ is the subset of
$\Gamma_1^{-}$ or, equivalently, if there exist two positive
constants $ N$ and $T$ such that $\|L^{n}\| \leq Ne^{-T n}$ for
all $n=0,1,2\dots$,

\item[(ii)]  \textit{expansive} if $\sigma(L)$ is the subset of
$\Gamma_1^{+}$ and

\item[(iii)] \textit{dichotomic} if $\sigma(L)$ does not intersect
the set $\Gamma_{1}$.
\end{itemize}
\vspace{.5cm}
\begin{remark}\label{re1}
If $L$ is a dichotomic matrix then there exists
$\eta \in {\{1,2,\dots ,\xi\}}$  such that
$$
|\lambda_1| \leq |\lambda_2| \leq \dots \leq |\lambda_{\eta}| < 1
< |\lambda_{\eta + 1}| \leq  \dots \leq |\lambda_\xi|.
$$
Having in mind the decomposition of $\mathbb{C}^m$ given by $(3.1)$
let us consider
$$
X_{1}= Y_1\oplus Y_2 \oplus \dots \oplus Y_{\eta}\quad \hbox{and}\quad
X_2= Y_{\eta + 1}\oplus Y_{\eta + 2} \oplus \dots \oplus
Y_\xi.
$$
Then  $\mathbb{C}^m = X_{1} \oplus  X_2. $
\end{remark}

Recall that a linear map $P: \mathbb{C}^{m}\to \mathbb{C}^{m}$ is called projection if $P^{2}=P.$ In the following theorem we give our first result.
\begin{theorem} Let $q>0$. If the matrix $L:=U(q,0)$ is dichotomic and there exists a projection $P$ commuting with $L$, $\Phi_{\mu}(q)$ and $\Psi_{\mu}(q)$ then for each $\mu \in \mathbb{R}$ and each non-zero vector $b \in \mathbb{C}^m$ the solutions of the following Cauchy problems
$$
\left\{
  \begin{array}{ll}
    \dot{X}(t)=A(t)X(t)+e^{i\mu t}Pb, \quad t\ge 0\\
X(0)=0,
  \end{array}
\right.\eqno{(A(t), \mu, Pb, 0)}$$
and
$$
\left\{
  \begin{array}{ll}
    \dot{X}(t)=-X(t)A(t)+e^{i\mu t}(I-P)b, \quad t\ge 0\\
X(0)=0,
  \end{array}
\right.\eqno{(-A(t), \mu, (I-P)b, 0)}$$
are bounded.
\end{theorem}
\begin{proof}
Assume that $L$ is dichotomic, then by Remark \ref{re1} we have a decomposition of $\mathbb{C}^m$, i.e.
$\mathbb{C}^m = X_{1} \oplus  X_2$. \\
We define $P : \mathbb{C}^m \to \mathbb{C}^m $  by $Px = x_{1}$,
where $x =  x_1 +  x_2$, such that $x_{1} \in X_{1}$ and $x_2 \in X_2$.
It is clear that $P$ is a projection. \\
Moreover for all $x \in \mathbb{C}^m $ and all $k\in\mathbb{Z_+}$, this yields
$$
PL^{k}x = P(L^{k}(x_{1} + x_2)) = P(L^{k}(x_{1}) +
L^{k}(x_2)) = L^{k}(x_{1}) = L^{k}Px.
$$
Hence
$PL^{k}=L^{k}P\,\,\,\text{  for all   } k\in\mathbb{Z_+}$. Also we have
$
P\Phi_{\mu}(q)x = P(\Phi_{\mu}(q)(x_{1} + x_2)) = P(\Phi_{\mu}(q)(x_{1}) +
\Phi_{\mu}(q)(x_2)) = \Phi_{\mu}(q)(x_{1}) = \Phi_{\mu}(q)Px
$ and similarly we conclude that $P\Psi_{\mu}(q)=\Psi_{\mu}(q)P.$ Now the solution of the
Cauchy problem
$(A(t), \mu, Pb, 0)$ is given by
$$\Phi_{(\mu,P,b)}(t)=\int_{0}^{t} U(t,s)e^{i\mu s}Pbds.$$

Let $n$ be the integer part of $\frac{t}{q}$ and let $r:=(t-qn)\in [0, q).$ Then
\begin{eqnarray*}
&&\int_{0}^{t}U(t,s)e^{i \mu s}Pb ds=\int_{0}^{qn +r}U(t,s)e^{i \mu s}Pb ds\\
&=&\int_{0}^{qn}U(t,s)e^{i \mu s}Pb ds+\int_{qn}^{qn +r}U(t,s)e^{i \mu s}Pb ds\\
&=&\int_{qn}^{qn +r}U(t,s)e^{i \mu s}Pb ds+\sum\limits_{k=0}^{n-1}\int_{qk}^{q(k+1)}U(qn +r,s)e^{i \mu s}Pb ds\\
&=&\int_{qn}^{qn +r}U(t,s)e^{i \mu s}Pb ds + U(r,0)\sum\limits_{k=0}^{n-1} \int_{qk}^{q(k+1)}U(qn,s)e^{i \mu s}Pb ds\\
&=&\int_{qn}^{qn +r}U(t,s)e^{i \mu s}Pb ds + U(r,0)\sum\limits_{k=0}^{n-1} \int_{0}^{q}U(qn ,qk+\tau)e^{i \mu (qk+ \tau)}Pb d\tau\\
&=&\int_{qn}^{qn +r}U(t,s)e^{i \mu s}Pb ds + U(r,0)\sum\limits_{k=0}^{n-1}e^{qi \mu k } \int_{0}^{q}U(q(n -k),\tau)e^{i \mu \tau}Pb d\tau\\
&=&\int_{qn}^{qn +r}U(t,s)e^{i \mu s}Pb ds+ U(r,0)\sum\limits_{k=0}^{n-1}e^{i \mu qk }U(q,0)^{n -k-1}\int_{0}^{q}U(q,\tau)e^{i \mu \tau}Pb d\tau\\
&=&I_{1}+I_{2}.
\end{eqnarray*}
where
$I_{1}=\int_{qn}^{qn +r}U(t,s)e^{i \mu s}Pb ds,$ and
$I_{2}=U(r,0)\sum\limits_{k=0}^{n-1}e^{i \mu qk }L^{n -k-1}\Phi_{\mu}(q)Pb$.\\
Now the family $\mathcal{U}$ has a growth bound and $0\leq t-s\leq r< q$, so we have
\begin{eqnarray*}
\|I_{1}\|&=&\left\|\int_{qn}^{qn +r}U(t,s)e^{i \mu s}Pb\,\,ds\right\|\\
&\leq& M \int_{qn }^{qn +r}e^{\omega(t-s)}\|Pb\|\\
&\leq& rM e^{q\omega}\|Pb\|\\
&\leq& qM e^{q\omega}\|Pb\|,
\end{eqnarray*}
where $\omega$ is a real number and $M\geq 1$. Hence $I_1$ is bounded.\\
Next let $ z_\mu= e^{i \mu  q}$, and $\Phi_{\mu}(q)b=l \in \mathbb{C}^{m}$ then
\[I_{2}=U(r,0)\big( L^{n-1}z_{\mu}^{0} + L^{n-2}z_{\mu}^{1} + \dots +
   L^{0}z_{\mu}^{n-1} \big)Pl\]
By our assumption we know that $L$ is dichotomic and $|z_\mu|= 1$ thus $z_\mu$ is contained in the resolvent set
of $L$ therefore the matrix $(z_\mu I - L )$  is an invertible matrix. Hence
\[I_{2}=U(r,0)(z_\mu I - L )^{-1} (z_\mu^{n} I - L^{n} )Pl.\]

Taking norm of both sides
\begin{eqnarray*}
\| I_{2}\| & \leq & \|U(r,0)(z_\mu I - L )^{-1}z_\mu^{n}Pl \| + \| U(r,0)(z_\mu I - L )^{-1}PL^{n}l\|\\
   &=& \|U(r,0) \|\|(z_\mu I - L )^{-1}\| \|Pl \| + \|U(r,0) \|\|(z_\mu I - L )^{-1}\|\| PL^{n}l\|.   \\
\end{eqnarray*}
 Using Theorem \ref{2.100}, we have
$$L^{n}l=\lambda_1^{n}p_1(n) + \lambda_2^{n}p_2(n) + \dots + \lambda_{\xi}^{n}p_{\xi}(n),$$
thus
$$PL^{n}l=\lambda_1^{n}p_1(n) + \lambda_2^{n}p_2(n) + \dots + \lambda_{\eta}^{n}p_{\eta}(n),$$
where each $p_i(n)$ are $\mathbb{C}^{m}$-valued polynomials with degree at most $(m_i -1)$ for any
$i \in \{ 1,2,\dots, \xi \}$.
From hypothesis we know that $ |\lambda_i| < 1$ for each $i \in \{ 1,2,\dots, \eta \}$. So
$\| PL^{n}l\| \rightarrow 0$ when $ n \rightarrow \infty$.
Thus $I_2$ is bounded, hence the solution of $(A(t), \mu, Pb, 0)$ is bounded.\\

Next, since the solution of the Cauchy problem
$(-A(t), \mu, (I-P)b, 0)$ is given by
$$\Psi_{(\mu,I-P,b)}(t)=\int_{0}^{t} V(t,s)e^{i\mu s}(I-P)b\,ds.$$
By similar method we obtain that
$$\Psi_{(\mu,I-P,b)}(t)=J_{1}+J_{2}$$
where $J_{1}=\int_{qn}^{qn +r}V(t,s)e^{i \mu s}(I-P)b ds$ and
$$J_{2}=V(r,0)(z_\mu^{0}L^{-(n-1)}+z_\mu^{1}L^{-(n -2)}+\dots+z_\mu^{n-1}L^{0})\Psi_{\mu}(q)(I-P)b.$$
Proceeding as before we can show that $J_{1}$ is bounded. Now for $J_{2}$ we have
since $PL=LP,$ therefore $(I-P)L=L(I-P).$ By our assumption we know that $L$ is invertible and since $L^{-1}$ is also dichotomic hence using the same arguments as above we have
\begin{eqnarray*}
J_{2}&=&V(r,0)(z_\mu I - L^{-1} )^{-1} (z_\mu^{n} I - L^{-n} ) \Psi_{\mu}(q)(I-P)b\\
 &=&V(r,0)(z_\mu I - L^{-1} )^{-1} (z_\mu^{n} I - L^{-n} )(I-P)\Psi_{\mu}(q)b.
\end{eqnarray*}
Taking norm of both sides we get
\begin{eqnarray*}
\| J_{2}\| &\leq& \|V(r,0) \|\|(z_\mu I - L^{-1} )^{-1}\| \|(I-P)\Psi_{\mu}(q)b \| \\ && + \|V(r,0) \|\|(z_\mu I - L^{-1} )^{-1}\|\| L^{-n}(I-P)\Psi_{\mu}(q)b\|.
\end{eqnarray*}
First we prove that $L^{-n}x \rightarrow 0$ as $n \rightarrow \infty$ for any $x\in X_{2}$. Since $(I-P)\Psi_{\mu}(q)b \in X_{2}$ the assertion would follows.
Now since $X_2= Y_{\eta + 1}\oplus Y_{\eta + 2} \oplus \dots \oplus Y_\xi.$
So any $x\in X_2$ can be written as a sum of $\xi-\eta $ vectors $y_{\eta+1}$, $y_{\eta+2}$, $\dots$ $y_{\xi}$. It would be sufficient to prove that
$L^{-n}y_i \rightarrow 0$ as $n \rightarrow \infty$ for any $i \in \{\eta+1, \eta+2,\dots,\xi \}$.
Let $Y\in \{Y_{\eta + 1},Y_{\eta + 2},\dots,Y_{\xi}\}$ say $Y=ker(L-\lambda I)^{\rho}$, where $\rho \geq 1$ is an integer number and $|\lambda|> 1$.
Consider $d_1 \in Y\backslash \{0\}$ such that $(L-\lambda I)d_1=0$ and let $d_2, d_3,\dots,d_\rho$ given by $(L-\lambda I)d_i=d_{i-1}$.
Then $A:=\{d_1,d_2,\dots, d_\rho\}$ is a basis in $Y$. So it is sufficient to prove that
$L^{-n}d_i \rightarrow 0$ as $n \rightarrow \infty$ for any $i \in \{1, 2,\dots,\rho \}$.
For $i=1,$ we have that
$L^{-n}d_1=\frac{1}{\lambda^{n}}d_1 \rightarrow 0$ as $n \rightarrow \infty$.
For $i=2,3,\dots,\rho$, denote $B_{n}=L^{-n}d_{i}$. Then $(L-\lambda I)^{\rho}B_{n}=0$, i.e.
\begin{equation*}
B_n - C_{\rho}^{1}B_{n-1}\alpha + C_{\rho}^{2}B_{n-2}\alpha^{2} +\dots +
C_{\rho}^{\rho}B_{n-\rho}\alpha^{\rho}=0 ,  \hbox{  for all  } n \geq \rho  \eqno(3.2)
\end{equation*}
where  $\alpha=\frac{1}{\lambda}.$\\
Passing for instance at the components, it follows that there exists a $\mathbb{C}^{m}$-valued polynomial
$P_\rho$ having degree at most $\rho-1$ and
verifying (3.2) such that $B_n= \alpha^{n}P_\rho(n).$ Thus $B_n\rightarrow0,\ when \ n \rightarrow \infty$ i.e.
$L^{-n}d_i\rightarrow 0$ for any  $i\in {\{1,2,\dots,\rho \}}.$ Thus $J_2$  is bounded.
\end{proof}

The converse statement of the above theorem is not straight forward and we need to put an extra condition i.e. the matrices $\Phi_{\mu}(q)$ and $\Psi_{\mu}(q)$ are invertible, at the end of the paper we have given an example which shows that the invertibility conditions on matrices $\Phi_{\mu}(q)$ and $\Psi_{\mu}(q)$ can not be removed. Due to this reason we put the converse statement of the above theorem as a new theorem which is stated as.
\begin{theorem}\label{3.3}
  Let there exists a projection $P$ commuting with $L$, $\Phi_{\mu}(q)$ and $\Psi_{\mu}(q)$ and let for each $\mu \in \mathbb{R}$ the matrices $\Phi_{\mu}(q)$ and $\Psi_{\mu}(q)$ are invertible then if for each real number $\mu$ and each non-zero vector $b\in \mathbb{C}^{m},$ the solutions of the Cauchy Problems $(A(t),\mu,Pb,0)$ and $(-A(t), \mu, (I-P)b, 0)$ are bounded then the Poincare map $L$ is dichotomic.
\end{theorem}
\begin{proof}
Suppose on contrary that the matrix $L$ is not dichotomic then
 $\sigma(L) \cap \Gamma_1 \neq \phi$. Let $ \omega \in \sigma(L)\cap \Gamma_1$ then there exists a non zero $y \in \mathbb{C}^{m}$ such that $Ly=\omega y$, it is easy to see that $L^k y= w^ky.$
 Here we have two cases:\\
Case 1: If $Py\neq0.$
Choose $\mu_1 \in \mathbb{R}$
such that $\omega=e^{i\mu_{1}q}$, then $L^k y= e^{i\mu qk}y.$
Since $\Phi_{\mu_{1}}(q)$ is invertible so there exists $b_{1}\in\mathbb{C}^m$ such that $\Phi_{\mu_{1}}(q)b_{1}=y.$ Then
\begin{eqnarray*}
\Phi_{(\mu_{1},P,b_{1})}(t)&=&\int_{qn}^{qn +r}U(t,s)e^{i \mu_{1} s}Pb_{1}\, ds+ U(r,0)\sum\limits_{k=0}^{n-1}e^{i\mu_{1}qk}PL^{n -k-1}y \\
&=&\int_{qn}^{qn +r}U(t,s)e^{i \mu_{1} s}Pb_{1}\, ds+ U(r,0)\sum\limits_{k=0}^{n-1}e^{i\mu_{1}qk}Pe^{i\mu_{1}q(n -k-1)}y \\
&=&\int_{qn}^{qn +r}U(t,s)e^{i \mu_{1} s}Pb_{1}\, ds+ U(r,0)\sum\limits_{k=0}^{n-1}e^{i\mu_{1}q(n-1)}Py \\
&=&\int_{qn}^{qn +r}U(t,s)e^{i \mu_{1} s}Pb_{1}\, ds+ U(r,0)ne^{i\mu_{1}q(n-1)}Py.
\end{eqnarray*}
Now clearly $U(r,0)ne^{i\mu_{1}q(n-1)}Py \to \infty \hbox{ as } n\to\infty$.
Hence there exist $\mu_{1}\in \mathbb{R}$ and $b_{1}\in \mathbb{C}^m$ such that $\Phi_{(\mu_{1},P,b_{1})}$ is unbounded.
Therefore contradiction arises.

Case 2: If $Py=0$ then surely $(I-P)y\neq0.$ Since $PL=LP$ therefore $(I-P)L=L(I-P).$
Choose $\mu_2 \in \mathbb{R}$
such that $\omega=e^{-i\mu_{2}q}$. In this case we note that $L^{-k}y=e^{i\mu_{2}qk}y$.
Also $\Psi_{\mu_{2}}(q)$ is invertible so there exists $b_{2}\in\mathbb{C}^m$ such that $\Psi_{\mu_{2}}(q)b_{2}=y.$
Now consider the solution of $(-A(t),\mu_2,b_2,0)$ we have
$$\Psi_{(\mu_{2},I-P,b_{2})}(t)=J_{1,\mu_2}+J_{2,\mu_2},$$
where
$$J_{1,\mu_2}=\int_{qn}^{qn +r}V(t,s)e^{i \mu_{2} s}(I-P)b_{2}\, ds,$$
and
\begin{eqnarray*}
J_{2,\mu_2}&=&V(r,0)\sum\limits_{k=0}^{n-1}e^{i\mu_{2}qk}L^{-(n -k-1)}\Psi_{\mu_{2}}(q)(I-P)b_{2}\\
&=&V(r,0)\sum\limits_{k=0}^{n-1}e^{i\mu_{2}qk}(I-P)L^{-(n -k-1)}y\\
&=&V(r,0)\sum\limits_{k=0}^{n-1}e^{i\mu_{2}qk}(I-P)e^{i\mu_2q(n -k-1)}y\\
&=&V(r,0)\sum\limits_{k=0}^{n-1}e^{i\mu_{2}q(n-1)}(I-P)y\\
&=&V(r,0)ne^{i\mu_{2}q(n-1)}(I-P)y.
\end{eqnarray*}
Clearly we see that
$J_{2,\mu_2}=V(r,0)n z_{\mu_{2}}^{n-1}(I-P)y\to \infty \text{ as } n\to\infty$.
Hence there exist $\mu_{2}\in \mathbb{R}$ and $b_{2}\in \mathbb{C}^m$ such that $\Psi_{(\mu_{2},I-P,b_{2})}(t)$ is unbounded.
Which is again an absurd. This completes the proof.
\end{proof}
The following theorem is taken from \cite{[SB]} which we used to obtained theorem 3.5.
\begin{theorem}\label{ues}
The matrix $L$ is stable if and only if for each $b\in \mathbb{C}^{m}$, the solution of
$(A(t),\mu,Pb,0)$ is bounded on $\mathbb{R}_{+}$ uniformly with respect to the parameter $\mu \in \mathbb{R}$, i.e
$$\sup_{\mu \in \mathbb{R}}\sup_{t\geq 0}\|\int_{0}^{t}U(t,s)e^{i\mu s}bds\|:=K(b)< \infty.$$
\end{theorem}
\begin{theorem}
The matrix $L$ is dichotomic if and only if there exists a projection $P$ such that for each vector $b\in \mathbb{C}^{m}$, the solutions of the Cauchy Problems $(A(t),\mu,Pb,0)$ and $(-A(t), \mu, (I-P)b, 0)$ are uniformly bounded on $\mathbb{R}_{+}$ with respect to the parameter $\mu \in \mathbb{R}$, i.e
$$\sup_{\mu \in \mathbb{R}}\sup_{t\geq 0}\|\int_{0}^{t}U(t,s)e^{i\mu s}Pbds\|:=K_{P}(b)< \infty,\eqno{(3.3)}$$ and
$$\sup_{\mu \in \mathbb{R}}\sup_{t\geq 0}\|\int_{0}^{t}V(t,s)e^{i\mu s}(I-P)bds\|:=K_{I-P}(b)< \infty.\eqno{(3.4)}$$
\end{theorem}
\begin{proof}
Suppose the matrix $L$ is dichotomic and let $L_{1}$ and $L_{2}$ be the restrictions of $L$ on $X_{1}$ and $X_{2}$ respectively. Consider the spectral decomposition of $\mathbb{C}^{m}$ as given in Remark \ref{re1}, that is we can write
$$\mathbb{C}^{m}=X_{1}\oplus X_{2}.$$
Then $L_{1}$ is stable on $X_{1}$  and $L_{2}^{-1}$ is stable on $X_{2}$.
Define the projection $P: \mathbb{C}^{m}\to \mathbb{C}^{m}$ as
$Px=x_{1}$ where $x=x_{1}+x_{2}$ such that $x_{1}\in X_{1}$ and $x_{2}\in X_{2}$. Then clearly
$P\mathbb{C}^{m}=X_{1}$ and $(I-P)\mathbb{C}^{m}=X_{2}$. \\
Since $Pb\in X_{1}$ for each $b\in \mathbb{C}^{m},$ therefore Theorem \ref{ues} implies that
$$\sup_{\mu \in \mathbb{R}}\sup_{t\geq 0}\|\int_{0}^{t}U(t,s)e^{i\mu s}Pbds\|:=K_{P}(b)< \infty.$$
Also $(I-P)b\in X_{2}$ for each $b\in \mathbb{C}^{m}$ then again Theorem \ref{ues} implies that
$$\sup_{\mu \in \mathbb{R}}\sup_{t\geq 0}\|\int_{0}^{t}V(t,s)e^{i\mu s}(I-P)bds\|:=K_{I-P}(b)< \infty.$$

Conversely let $P$ be the projection for which (3.3) and (3.4) are satisfied. Assume that $P\mathbb{C}^{m}=W_{1}$ and $(I-P)\mathbb{C}^{m}=W_{2}$. Then clearly $\mathbb{C}^{m}=W_{1}\oplus W_{2}$.
So by (3.3) and using Theorem \ref{ues} we have $L$ is stable on $W_{1}$.
Similarly by (3.4) and again using Theorem \ref{ues} we obtain that
$L^{-1}$ is stable on $W_{2}$.
Hence $L$ is dichotomic on $\mathbb{C}^{m}$.
\end{proof}

Now we will present an example which shows that in Theorem \ref{3.3} the invertibility condition on the matrices $\Phi_{\mu}(q)$ and $\Psi_{\mu}(q)$
can not be removed.
\begin{example} Let
$$\Phi(t)=\left(
    \begin{array}{cc}
      \cos t & \sin t \\
      -\sin t & \cos t \\
    \end{array}
  \right)
,$$
then
$$\Phi^{-1}(t)=\Phi(-t)=\left(
    \begin{array}{cc}
      \cos t & -\sin t \\
      \sin t & \cos t \\
    \end{array}
  \right).$$
  So in this case the evolution family ${\mathcal{U}}=\{U(t,s), t,s\in\mathbb{R}_{+}\},$ is given by
  $$U(t,s)=\Phi(t)\Phi^{-1}(s)=\left(
    \begin{array}{cc}
      \cos(t-s) & \sin(t-s) \\
    -\sin(t-s) & \cos(t-s) \\
    \end{array}
  \right).$$
  Since $\sin t$ and $\cos t$ are $2\pi$-periodic functions so this evolution family is $2\pi$-periodic.
Now
$$U(2\pi,s)=\left(
    \begin{array}{cc}
      \cos s & -\sin s \\
      \sin s & \cos s  \\
    \end{array}
  \right)
,$$
and
$$\Phi_{\mu}(2\pi)=\int_{0}^{2\pi}U(2\pi,s)e^{i \mu s}ds.$$
Choose $\mu=0,$ we have
$$\Phi_{0}(2\pi)=\int_{0}^{2\pi}U(2\pi,s)ds=\left(
                                              \begin{array}{cc}
                                                0 & 0 \\
                                                0 & 0 \\
                                              \end{array}
                                            \right)
$$
which is not invertible.
The solution of the Cauchy problem $(A(t),\mu,Pb,0)$ is given by
$$\Phi_{\mu}(t)Pb=\int_{2\pi n}^{2\pi n +r}U(t,s)e^{i \mu s}Pb ds+ U(r,0)\sum\limits_{k=0}^{n-1}z_{\mu}^{k}U(2\pi,0)^{n -k-1}\Phi_{\mu}(2\pi)Pb,$$
where $r\in [0,2\pi).$
Now the family has growth bound and $0\leq t-s<2 \pi$, so we have
\begin{eqnarray*}
\left\|\int_{2\pi n}^{2\pi n +r}U(t,s)e^{i \mu s}Pb ds\right\|&\leq& rMe^{2\pi \omega } \|Pb\|\\
 &\leq& 2\pi Me^{2\pi \omega} \|Pb\|<\infty,
\end{eqnarray*}
where $\omega \in \mathbb{R}$ and $M\geq 1$.\\
Next we have since

$$U(2\pi,0)=\left(
              \begin{array}{cc}
                1 & 0 \\
                0 & 1 \\
              \end{array}
            \right)
,$$
therefore
$$\sum\limits_{k=0}^{n-1}z_{\mu}^{k}U(2\pi,0)^{n -k-1}=\left(
              \begin{array}{cc}
                \sum\limits_{k=0}^{n-1}z_{\mu}^{k} & 0 \\
                0 & \sum\limits_{k=0}^{n-1}z_{\mu}^{k} \\
              \end{array}
            \right)
.$$
For $z_{\mu}\neq1,$ we obtain
$$
| \sum\limits_{k=0}^{n-1}z_{\mu}^{k}| = \frac{\big|
z_{\mu }^{k}-1\big|}{\big|1 - z_{\mu }\big|} \leq \frac{2}{| 1 - z_{\mu }|}.
$$
So for the corresponding values of  $\mu \in \mathbb{R}$ and each $b\in \mathbb{C}^{2}$ the solution is bounded. If $z_\mu = 1$ i.e. $\mu=0,$ then

$$\begin{array}{rcl} \sum\limits_{k=0}^{n-1}z_{\mu}^{k}U(2\pi,0)^{n -k-1}\Phi_{0}(2\pi)
              &=& \left(
                   \begin{array}{cc}
                     n & 0  \\
                       0 &  n \\
                   \end{array}
                 \right)
                 \left(
                   \begin{array}{cc}
                     0 & 0  \\
                       0 & 0  \\
                   \end{array}
                 \right)=
                 \left(
                   \begin{array}{cc}
                      0& 0  \\
                       0 & 0 \\
                   \end{array}
                 \right).
                 \end{array}$$
Thus the solution is bounded.\\
Also we have
$$V(t,s)=U^{-1}(t,s)=\left(
    \begin{array}{cc}
      \cos(t-s) & -\sin(t-s) \\
    \sin(t-s) & \cos(t-s) \\
    \end{array}
  \right).$$
Similarly as above we can see that $\Psi_{0}(2\pi)$ is not invertible and the solution $\Psi_{(\mu,I-P,b)}(t)$ is bounded for each $\mu\in \mathbb{R}$ and $b\in\mathbb{C}^{2}.$
But $1\in\sigma(U(2\pi, 0)),$ i.e. the matrix $U(2\pi, 0)$ is not dichotomic.
\end{example}


\begin{thebibliography}{9}

\vspace{0.5cm}

\bibitem{[SB]} S. Arshad, C. Bu\c se and Olivia Saierli, {\it Connections between exponential stability and boundedness of solutions of a couple of differential time depending and periodic systems}, Electronic Journal of Qualitative Theory of Differential Equations,  (2011), No. 90, pp. 1-16.

\bibitem{[BZ]}  C. Buse,  Akbar Zada, Dichotomy and bounded-ness of solutions for some discrete Cauchy problems, Proceedings of IWOTA- 2008, Operator Theory, Advances and Applications, (OT) Series Birkh\"auser Verlag, Eds: J.A. Ball, V. Bolotnikov, W. Helton, L. Rodman, T. Spitkovsky, Vol. {\bf 203}, 165-174, (2010).


\bibitem{[ABNZ11]} S. Arshad, C. Bu\c se, A. Nosheen and A. Zada, {\it Connections between the stability of a Poincare map and boundedness of certain associate sequences}, Electronic Journal of Qualitative Theory of Differential Equations,  (2011), No. 16, pp. 1-12.


\bibitem{[CY]} C. Chicone, Y.Latushkin,  "Evolution Semigroups in Dynamical Systems and Differential Equations", Amer. Math. Soc., Math. Surv. and Monographs, No. 70 (1999).



\bibitem {[BZ09]} C. Bu\c se and A. Zada, {\it Boundedness and exponential stability for periodic time dependent systems},  Electronic Journal of Qualitative Theory of Differential Equations, {\bf 37}, pp 1-9 (2009).



\bibitem{[Z1]} Akbar Zada, A characterization of dichotomy in terms of boundedness of solutions for some Cauchy problems, Electronic Journal of Differential Equations, Vol.(2008), No. 94, 1-5.


\bibitem{[BCDS]} C. Bu\c se, P. Cerone, S. S. Dragomir and A. Sofo, {\it Uniform stability of periodic discrete system in Banach spaces,} J. Difference Equ. Appl. \textbf{11}, No .12 (2005) 1081-1088.


\end{thebibliography}
\end{document}